\documentclass[twoside,a4paper,11pt,reqno,lref] {amsart}
\usepackage{amsmath,mathrsfs}
\numberwithin{equation}{section}
\renewcommand{\theequation}{\thesection.\arabic{equation}}

\begin{document}
\title[Volume Growth and Curvature Decay]
{Volume Growth and Curvature Decay of Complete Positively Curved K\"{a}hler Manifolds}

%    Information for first author
\author{Xiaoyong Fu}
\address{Department of Mathematics, Zhongshan University,
Guangzhou 510275, China
}
\email{mcsfxy@mail.sysu.edu.cn}
%    \thanks will become a 1st page footnote.
%\thanks{}

%    Information for second author
\author{Zhenglu Jiang}
%    Address of record for the research reported here
\address{Department of Mathematics, Zhongshan University,
Guangzhou 510275, China}
\email{mcsjzl@mail.sysu.edu.cn}
%\thanks{}
%    Information for third author
%\thanks{}

%\addtocounter{footnote}{0}
%\pagestyle{myheadings}
%\markboth{Z.~Jiang, X.~Fu}{Volume Growth and Curvature Decay}

%\noindent{\footnotesize Submitted \today}
%\baselineskip 5pt
%\vskip0.6cm

\date{\today.}

%\dedicatory{This paper is dedicated to our advisors.}


\begin{abstract}
This paper constructs a class of complete K\"{a}hler metrics of positive holomorphic  
sectional curvature on ${\bf C}^n$ and finds that  the constructed metrics satisfy the following properties: As 
the geodesic distance $\rho\to\infty,$ 
 the volume of geodesic balls grows like $O(\rho^{\frac{2(\beta+1)n}{\beta+2}})$  
 and the Riemannian scalar curvature 
decays like $O(\rho^{-\frac{2(\beta+1)}{\beta+2}}),$ where $\beta\geq 0.$ 
\end{abstract}
\maketitle

\section{Introduction}
\label{intro} 
  We are concerned with the volume growth and curvature decay of complex $n$-dimensional ($n\geq 2$) complete
 noncompact K\"{a}hler manifolds (denoted by $M$) with positive holomorphic sectional curvature.  
For convenience, throughout this paper, $\rho=\rho(x_0,x)$ represents the geodesic distance 
from a fixed point $x_0$ to $x$ in $M,$   
 $V(B(x_0,\rho))$ denotes the volume of the geodesic ball $B(x_0,\rho)$ 
centered at $x_0$ with radius $\rho$  and 
$R(x)$  the Riemannian scalar curvature at $x.$ 
When $M$ is assumed to be of positive holomorphic bisectional curvature, 
it is known by the classical Bishop volume comparison theorem that 
$V(B(x_0,\rho)) \leq \omega(n)\rho^{2n}$  for any $x_0\in M,$  
 where $\omega(n)$ is the volume 
of the standard unit ball $B^{2n}$ in ${\bf R}^{2n}.$ Also, if $V(B(x_0,\rho)) \geq C\rho^{2n}$ for 
some positive constant $C,$ then the scalar curvature decays quadratically in the average sense. 
This is conjectured by Yau \cite{y}
and confirmed by Chen and Zhu \cite{cz}. 
On the other hand, it has been shown by Chen and Zhu \cite{cz} that \\ 
 \centerline{$V(B(x_0,\rho)) \geq C\rho^{n}$ for any $x_0\in M,$}  
 where $C$ is a positive constant 
depending only on $x_0$ and $M,$ and that the scalar curvature decays at least linearly  in the average sense.
In view of the above results, 
a complex $n$-dimensional K\"{a}hler manifold with positive holomorphic bisectional curvature is said to be of maximal volume growth 
if $V(B(x_0,\rho)) =O(\rho^{2n})$ and of minimal volume growth if $V(B(x_0,\rho)) =O(\rho^{n}).$

 There have been some examples of positively curved complete noncompact  K\"{a}hler manifold
 with  maximal (respectively, minimal) volume growth and quadratic (respectively, linear) curvature 
decay. In \cite{pfk}, Klembeck constructed a complete, rotationally symmetric K\"{a}hler metric $g$
of positive holomorphic sectional curvature on ${\bf C}^n$ with minimal volume growth and linear curvature decay. 
More precisely, the volume of the geodesic ball $B(0,\rho)$ with respect to 
the metric $g$ grows like 
$O(\rho^{n})$ and the scalar curvature $R(x)$ of the metric $g$ decays like $O(\rho^{-1}).$ 
Rotationally symmetric, complete gradient  K\"{a}hler-Ricci solitons of positive 
 Riemannian sectional curvature on ${\bf C}^n$ 
have been found by Cao (\cite{c1}, \cite{c2}) and it has been known that
such solitons can be divided into two branches: one is of maximal volume growth and quadratic curvature 
decay, and the other of minimal volume growth and linear curvature decay.

By the above analysis, it is natural to ask whether there exists a complete noncompact K\"{a}hler
manifold of positive holomorphic sectional curvature  satisfying \\
 \centerline{$V(B(x_0,\rho)) = O(\rho^{n(1+\epsilon )})$ and $R(x) = O(\rho^{-(1+\epsilon )})$ 
for any fixed  $\epsilon\in (0,1).$} 
To our best knowledge, examples of such manifolds seem to be elusive from the literature.
In this paper, we shall construct a class of complete K\"{a}hler metrics $g$  of positive holomorphic  
sectional curvature on ${\bf C}^n$ and find that  
the constructed metrics satisfy the following properties: As 
the geodesic distance $\rho\to\infty,$ 
 the volume of geodesic balls grows like $O(\rho^{\frac{2(\beta+1)n}{\beta+2}})$  
 and the Riemannian scalar curvature 
decays like $O(\rho^{-\frac{2(\beta+1)}{\beta+2}}),$ where $\beta\geq 0.$ 

   Finally, we should mention the work of Ni et al.~\cite{nst} for some results on the relation
 between volume growth and curvature decay of complete noncompact K\"{a}hler manifold
with positive holomorphic bisectional curvature , and refer to the book of Greene and Wu \cite{gw} 
and  Kobayashi and Nomizu \cite{kn} for background on 
 differential geometry.
  
\section{Method of Construction}
In this section, we shall give the conditions under which there exists 
a complete rotationally symmetric K\"{a}hler metric on ${\bf C}^n$ with  positive  
 holomorphic  sectional curvature. 

Recall that a complex $n$-dimensional complete K\"{a}hler manifold
$M$ is of positive holomorphic sectional curvature if $\sum R_{i\bar{j}k\bar{l}}a_i\bar{a}_ja_k\bar{a}_l>0$ 
for all nonzero $n$-tuples $(a_1,a_2,\cdots,a_n)$ of complex numbers, 
where $R_{i\bar{j}k\bar{l}}$ is the components of the curvature tensor. 

Let $r^2=\sum\limits_{i=1}^nz_i\bar{z_i}$ on ${\bf C}^n.$ 
As in \cite{pfk}, we consider only rotationally symmetric metrics 
$g_{i\bar{j}}=\frac{\partial^2f(r^2)}{\partial z_i\partial\bar{z_j}}$ on ${\bf C}^n,$ 
derived from a global potential function, 
where $r\to f(r^2)\in C^\infty({\bf R}).$  Clearly, 
\begin{equation}
g=f^\prime(r^2)\left(\sum\limits_{i=1}^ndz_id\bar{z_i}\right)+
f^{\prime\prime}(r^2)\left(\sum\limits_{i=1}^n\bar{z_i}dz_i\right)\left(\sum\limits_{i=1}^nz_id\bar{z_i}\right).
\label{metric}
\end{equation}
Thus $g$ is a complete metric on ${\bf C}^n$ if $f$ satisfies the following two conditions: 
\begin{list}{}{}
\item[(i)] $f^\prime(r^2)>0$ and $f^\prime(r^2)+r^2f^{\prime\prime}(r^2)>0$ for all $r,$
\item[(ii)] $\int_0^\infty\sqrt{f^\prime(r^2)+r^2f^{\prime\prime}(r^2)}dr$ diverges.  
\end{list}
Since $g$ is a rotationally symmetric K\"{a}hler metric, without loss of generality, we may 
restrict the computation of the curvature tensor $R_{j\bar{k}l\bar{m}}$ of $g$ to the complex line 
$$L=\{ z_i=0|i>1\},$$ 
thus obtaining  
\begin{eqnarray}
R_{j\bar{k}l\bar{m}}&= &-\left(\frac{\partial^2 g_{j\bar{k}}}{\partial z_l\partial \bar{z}_m}
-\sum\limits_{p,q}g^{\bar{p}q}\frac{\partial g_{j\bar{p}}}{\partial z_l}
\frac{\partial g_{\bar{k}q}}{\partial \bar{z}_m}\right)  \nonumber \\ 
&=&-f^{\prime\prime}(r^2)
(\delta _{j\bar{k}}\delta _{l\bar{m}}  
+ \delta _{j\bar{m}}\delta _{l\bar{k}}) \nonumber \\ 
& & { }-r^2[f^{\prime\prime\prime}(r^2)-\frac{(f^{\prime\prime}(r^2))^2}
{f^{\prime}(r^2)}](\delta _{j\bar{k}1}\delta _{l\bar{m}}
+\delta _{j\bar{m}1}\delta _{l\bar{k}}
+\delta _{l\bar{m}1}\delta _{j\bar{k}}
+\delta _{l\bar{k}1}\delta _{j\bar{m}}) \nonumber \\
& & { }   -[r^4f^{\prime\prime\prime\prime}(r^2)-r^2\frac{(2f^{\prime\prime}(r^2)
+r^2f^{\prime\prime\prime}(r^2))^2}{f^\prime(r^2)+r^2f^{\prime\prime}(r^2)} 
+4r^2\frac{(f^{\prime\prime}(r^2))^2}{f^{\prime}(r^2)}]\delta _{j\bar{k}l\bar{m}1}.    \label{r01} 
\end{eqnarray}
Let $(a_1,a_2,\cdots,a_n)$ be any complex $n$-tuples  of complex numbers.
Then 
$$\sum R_{i\bar{j}k\bar{l}}a_i\bar{a}_ja_k\bar{a}_l=-(2A+4B+C)|a_1|^4-4(A+B)|a_1|^2(\sum\limits_{j=2}^n|a_j|^2)-2A(\sum\limits_{j=2}^n|a_j|^2)^2$$
where $$A=f^{\prime\prime}(r^2),$$ 
$$B=r^2[f^{\prime\prime\prime}(r^2)-\frac{(f^{\prime\prime}(r^2))^2}
{f^{\prime}(r^2)}],$$ 
$$C=r^4f^{\prime\prime\prime\prime}(r^2)-r^2\frac{(2f^{\prime\prime}(r^2)
+r^2f^{\prime\prime\prime}(r^2))^2}{f^\prime(r^2)+r^2f^{\prime\prime}(r^2)} 
+4r^2\frac{(f^{\prime\prime}(r^2))^2}{f^{\prime}(r^2)}.$$
Thus, in addition to (i) and (ii),  the following conditions for $g_{i\bar{j}}=\frac{\partial^2f(r^2)}{\partial z_i\partial\bar{z_j}}$ 
to be a complete K\"{a}hler metric on ${\bf C}^n$ of strictly positive holomorphic sectional curvature are required: 
\begin{list}{}{}
\item[(iii)] $A<0,$ i.e., $f^{\prime\prime}(r^2)<0$ for all $r,$
\end{list} 
\begin{list}{}{}
\item[(iv)] $2A+4B+C<0,$ i.e., $\frac{1}{4r}\frac{\partial}{\partial r}\left(r\frac{\partial}{\partial r}\ln(f^\prime(r^2)+r^2f^{\prime\prime}(r^2))\right)< 0,$
\item[(v)] $A+B<0,$ i.e., $f^{\prime\prime}(r^2)+r^2f^{\prime\prime\prime}(r^2)-r^2\frac{(f^{\prime\prime}(r^2))^2}{f^{\prime}(r^2)}<0.$
\end{list}
Therefore a complete metric of positive  holomorphic sectional 
curvature on ${\bf C}^n$ can be generated from a function $f\in C^\infty({\bf R})$ satisfying (i)-(v). 
In the next section, we shall introduce a class of functions in $C^\infty({\bf R})$ satisfying (i)-(v). 

\section{Metrics of Positive Curvature}
In this section, we shall show a family of complete K\"{a}hler metrics of positive  
holomorphic sectional curvature on ${\bf C}^n$ and discuss the curvature decay and volume growth of 
these metrics. 

Let us consider the following function 
\begin{equation}
f(r^2)=\frac{1}{(\beta+1)\alpha^\beta}\int_0^{r^2}\left\{\frac{[\alpha+\ln(1+x)]^{\beta+1}-\alpha^{\beta+1}}{x}\right\}dx\in C^\infty({\bf R}) \label{ex02} 
\end{equation}
where $\alpha>\beta\geq 0.$ Then 
$$f^\prime(r^2)=\frac{1}{(\beta+1)\alpha^\beta}\left\{\frac{[\alpha+\ln(1+r^2)]^{\beta+1}-\alpha^{\beta+1}}{r^2}\right\},$$
$$f^{\prime\prime}(r^2)=\frac{-[\alpha+\ln(1+r^2)]^{\beta+1}}{(\beta+1)\alpha^\beta r^4}
+\frac{[\alpha+\ln(1+r^2)]^{\beta}}{\alpha^\beta r^2(1+r^2)}+\frac{\alpha}
{(\beta+1) r^4}$$
and 
$$f^\prime(r^2)+r^2f^{\prime\prime}(r^2)=\frac{[\alpha+\ln(1+r^2)]^\beta}{\alpha^\beta(1+r^2)}.$$
It follows that $f$ satisfies (i), (ii) and (iii) (see \ref{appA}).  Also, 
$$\frac{1}{4r}\frac{\partial}{\partial r}\left(r\frac{\partial}{\partial r}\ln(f^\prime(r^2)+r^2f^{\prime\prime}(r^2))\right)
$$$$=-\frac{\alpha(\alpha-\beta)+\beta r^2+(2\alpha-\beta)\ln(1+r^2)+\ln^2(1+r^2)}{(1+r^2)^2[\alpha+\ln(1+r^2)]^2}$$
and 
$$f^{\prime\prime}(r^2)+r^2f^{\prime\prime\prime}(r^2)-r^2\frac{(f^{\prime\prime}(r^2))^2}{f^{\prime}(r^2)} $$
\begin{equation}
=-\frac{\beta\alpha^{\beta+1}e^{y-\alpha}-\beta\alpha^{\beta+1}-y^{\beta+2}+y^{\beta+1}e^{y-\alpha}-y^{\beta+1}+\alpha^{\beta+1}y}
{y^{-\beta}r^2(1+r^2)^2\{[\alpha+\ln(1+r^2)]^{\beta+1}-\alpha^{\beta+1}\}} \label{con5proof}
\end{equation} 
which yield (iv) and (v) (see \ref{appA}). Here, $y=\alpha+\ln(1+r^2).$ 

Therefore $g_{i\bar{j}}=\frac{\partial^2f(r^2)}{\partial z_i\partial\bar{z_j}}$ with $f(r^2)$ 
defined by (\ref{ex02}) 
is a class of complete K\"{a}hler metrics of positive sectional curvatures on ${\bf C}^n.$ 

Now we turn to the computation of the volume growth and 
scalar curvature decay of ${\bf C}^n$ equipped with the metric 
$g_{i\bar{j}}=\frac{\partial^2f(r^2)}{\partial z_i\partial\bar{z_j}}$  with $f(r^2)$ defined by (\ref{ex02}). 
First, let us estimate the volume growth of ${\bf C}^n.$ 
As before, our computation is restricted on $L.$   
Inserting (\ref{ex02}) into (\ref{metric}), we have   
\begin{equation}
g=\left\{\frac{[\alpha+\ln(1+r^2)]^\beta}{\alpha^\beta(1+r^2)}\right\}dz_1d\bar{z_1}+
\left\{\frac{[\alpha+\ln(1+r^2)]^{\beta+1}-\alpha^{\beta+1}}{(\beta+1)\alpha^\beta r^2}\right\}
\left(\sum\limits_{i=2}^ndz_id\bar{z_i}\right).   
\label{m02} 
\end{equation} 
And the volume form $\omega^n$ of (\ref{m02}) is given by 
$$\omega^n=\left(\frac{\sqrt{-1}}{2}\right)^n\left\{\frac{[\alpha+\ln(1+r^2)]^\beta}{\alpha^\beta(1+r^2)}\right\} $$
$$\times\left\{\frac{[\alpha+\ln(1+r^2)]^{\beta+1}-\alpha^{\beta+1}}{(\beta+1)\alpha^\beta r^2}\right\}^{n-1}
dz_1\wedge d\bar{z_1}\wedge dz_2\wedge d\bar{z_2}\wedge\cdots \wedge dz_nd\wedge\bar{z_n}.$$ 
Since the geodesic distance function $\rho$ from the origin of ${\bf C}^n$ is given by 
$$\rho=\int_{0}^r\sqrt{\frac{[\alpha+\ln(1+t^2)]^\beta}{\alpha^\beta(1+t^2)}}dt$$ and satisfies 
\begin{equation}
\rho=O(\ln^{\frac{\beta+2}{2}}r)~~(r\to\infty),
\label{gdp}
\end{equation} 
the volume growth $V(B(0,\rho))$ of geodesic ball $B(0,\rho)$ of ${\bf C}^n$ 
equipped with (\ref{m02}) is   
$$V(B(0,\rho))=\int_{B_E(0,r)}\omega^n 
=\int_{S^{2n-1}(1)}\left[\int_{0}^{r}\left(\frac{\sqrt{-1}}{2}\right)^{n+1}
\left\{\frac{[\alpha+\ln(1+t^2)]^\beta}{\alpha^\beta(1+t^2)}\right\}\right.$$ $$ 
\left.\times\left\{\frac{[\alpha+\ln(1+t^2)]^{\beta+1}-\alpha^{\beta+1}}{(\beta+1)\alpha^\beta t^2}\right\}^{n-1}t^{2n-1}dt\right]d\theta$$
$$=O(\ln^{(\beta+1)n}r)~~~\hbox{as}~r\to\infty $$
$$=O(\rho^{\frac{2(\beta+1)n}{\beta+2}})~~~\hbox{as}~\rho\to\infty ,$$ 
where $B_E(0,r)$ is the Euclidean ball corresponding to the geodesic ball $B(0,\rho).$

To determine the rate of curvature decay, without loss of generality, we may 
restrict the computation of the scalar curvature $R=\sum\limits_{i,j}g^{i\bar{j}}R_{i\bar{j}}$ of $g$ to the complex line 
$L=\{ z_i=0|i>1\}.$ Then the Ricci curvature, denoted by $\hbox{Ric},$  is as follows:  
$$\hbox{Ric}=-\sqrt{-1}
\left\{\beta\partial\bar{\partial}\ln\left[\alpha+\ln(1+r^2)\right]-\partial\bar{\partial}\ln(1+r^2)\right\}$$
$$-\sqrt{-1}(n-1)\left\{\partial\bar{\partial}\ln\left([\alpha+\ln(1+r^2)]^{\beta+1}-\alpha^{\beta+1}\right)-\partial\bar{\partial}\ln r^2\right\}$$
$$=-\sqrt{-1}\left\{\frac{\beta}{\alpha+\ln(1+r^2)}-1\right\}
\frac{dz_1\wedge d\bar{z_1}+(1+r^2)(\sum\limits_{i=2}^ndz_i\wedge d\bar{z_i})}{(1+r^2)^2}$$
$$+\sqrt{-1}\frac{\beta r^2}{[\alpha+\ln(1+r^2)]^2(1+r^2)^2}dz_1\wedge d\bar{z_1}$$
$$-\sqrt{-1}(n-1)\frac{(\beta+1)[\alpha+\ln(1+r^2)]^\beta}{[\alpha+\ln(1+r^2)]^{\beta+1}-\alpha^{\beta+1}}
\frac{\sum\limits_{i=1}^ndz_i\wedge d\bar{z_i}}{1+r^2}$$
$$-\sqrt{-1}(n-1)\left\{\frac{(\beta+1)[\alpha+\ln(1+r^2)]^{\beta-1}[\beta-\alpha-\ln(1+r^2)]}{[\alpha+\ln(1+r^2)]^{\beta+1}-\alpha^{\beta+1}}\right.$$
$$\left. -\frac{(\beta+1)^2[\alpha+\ln(1+r^2)]^{2\beta}}{([\alpha+\ln(1+r^2)]^{\beta+1}-\alpha^{\beta+1})^2}\right\}
\frac{r^2dz_1\wedge d\bar{z_1}}{(1+r^2)^2}+\sqrt{-1}(n-1)\frac{\sum\limits_{i=2}^ndz_i\wedge d\bar{z_i}}{r^2}.$$
This implies that
$$R_{1\bar{1}}=\left\{\frac{\beta}{\alpha+\ln(1+r^2)}-1\right\}\frac{1}{(1+r^2)^2}$$
$$+\frac{(n-1)(\beta+1)[\alpha+\ln(1+r^2)]^\beta}{[\alpha+\ln(1+r^2)]^{\beta+1}-\alpha^{\beta+1}}
\frac{1}{1+r^2}-\frac{\beta r^2}{[\alpha+\ln(1+r^2)]^2(1+r^2)^2}$$
$$+\frac{(n-1)(\beta+1)[\alpha+\ln(1+r^2)]^{\beta-1}[\beta-\alpha-\ln(1+r^2)]}{[\alpha+\ln(1+r^2)]^{\beta+1}-\alpha^{\beta+1}}
\frac{r^2}{(1+r^2)^2}$$
$$-\frac{(n-1)(\beta+1)^2[\alpha+\ln(1+r^2)]^{2\beta}}{([\alpha+\ln(1+r^2)]^{\beta+1}-\alpha^{\beta+1})^2}\frac{r^2}{(1+r^2)^2}$$
and that for $i\geq 2,$
$$R_{i\bar{i}}=\left\{\frac{\beta}{\alpha+\ln(1+r^2)}-1\right\}\frac{1}{1+r^2}-\frac{n-1}{r^2}$$
$$+\frac{(n-1)(\beta+1)[\alpha+\ln(1+r^2)]^\beta}{[\alpha+\ln(1+r^2)]^{\beta+1}-\alpha^{\beta+1}}
\frac{1}{1+r^2}.$$ 
Also, $R_{i\bar{j}}=0$ for $i\not=j.$ 
As a consequence, the scalar curvature $R=\sum\limits_{i=1}^ng^{i\bar{i}}R_{i\bar{i}}$ of ${\bf C}^n$ equipped with the rotationally  
symmetric metrics $g_{i\bar{j}}=\partial^2f(r^2)/\partial z_i\partial\bar{z_j}$  
with $f(r^2)$ defined by (\ref{ex02}) decays like
\begin{equation}
R=O\left(\frac{1}{\ln^{\beta+1}r}\right) \hbox{ as }r\to\infty. 
\label{sk01}
\end{equation} 
By (\ref{gdp}),  (\ref{sk01}) can be written as 
\begin{equation}
R=O\left(\rho^{-\frac{2(\beta+1)}{\beta+2}}\right) \hbox{ as } \rho\to\infty. 
\label{sk02}
\end{equation}

%\numberwithin{equation}{section}
\setcounter{equation}{0}
\setcounter{section}{0}
\setcounter{figure}{0}
\renewcommand{\thefigure}{\Alph{section}.\arabic{figure}}
\renewcommand{\thesection}{{Appendix \Alph{section}}}
\renewcommand{\theequation}{\Alph{section}.\arabic{equation}}
%\appendix
%\numberwithin{equation}{section}
%\setcounter{equation}{0}
%\renewcommand{\theequation}{B.\arabic{equation}}
\section{}
\label{appA}
Put $G(x)=[\alpha+\ln(1+x)]^{\beta+1}(1+x)-(\beta+1)x[\alpha+\ln(1+x)]^\beta -\alpha^{\beta+1}(1+x).$ Then 
$f^{\prime\prime}(r^2)=-\frac{G(r^2)}{(\beta+1)\alpha^\beta r^4(1+r^2)}.$ 
Since $\lim\limits_{r^2\rightarrow 0}f^{\prime\prime}(r^2)=-\frac{\alpha-\beta}{2\alpha}<0$ for $\alpha>\beta\geq 0,$ 
(iii) holds if $G(x)> 0$ for all $x> 0.$ 
Since $G(0)=G^\prime(0)=0,$ it suffices to show that 
$G^{\prime\prime}(x)> 0$ when $\alpha>\beta\geq 0.$ Indeed, if $\alpha>\beta\geq 0$ and $x\geq 0,$ then 
$$G^{\prime\prime}(x)=\frac{[(2\alpha-\beta)+2\alpha x]\ln(1+x)}{(\beta+1)^{-1}(1+x)(\alpha+\ln(1+x))^{-\beta}} $$
$$+\frac{\alpha(\alpha-\beta)+[\alpha^2-\beta^2+\beta]x}{(\beta+1)^{-1}(1+x)(\alpha+\ln(1+x))^{-\beta}}  
+\frac{\ln^2(1+x)}
{(\beta+1)^{-1}(\alpha+\ln(1+x))^{-\beta}}> 0.$$

To prove (v), let $H(y)=\beta\alpha^{\beta+1}e^{y-\alpha}-\beta\alpha^{\beta+1}-y^{\beta+2}+y^{\beta+1}e^{y-\alpha}-y^{\beta+1}+\alpha^{\beta+1}y.$ 
Since 
$$\lim\limits_{r^2\to 0}\left[f^{\prime\prime}(r^2)
+r^2f^{\prime\prime\prime}(r^2)
-r^2\frac{(f^{\prime\prime}(r^2))^2}{f^{\prime}(r^2)}\right]=-\frac{\alpha-\beta}{2\alpha}<0,$$
 it is easily known from (\ref{con5proof}) that (v) holds if $H(y)> 0$ for all $y> \alpha.$ 
Since $H(\alpha)=H^\prime(\alpha)=0,$ it suffices to prove that 
 $$H^{\prime\prime}(y)=\beta\alpha^{\beta+1}e^{y-\alpha}+y^\beta[ye^{y-\alpha}-\beta(\beta+1)]
+y^{\beta-1}[e^{y-\alpha}-1]\beta(\beta+1)+2y^\beta [e^{y-\alpha}-1](\beta+1)> 0$$ 
for all $y> \alpha$  
 when $\alpha>\beta\geq 0.$ 
 
Let 
$$I(y)=\beta\alpha^{\beta+1}e^{y-\alpha}+y^\beta[ye^{y-\alpha}-\beta(\beta+1)].$$
Then we need to show that $I(y)> 0$ for all $y\geq \alpha$  
 when $\alpha>\beta\geq 0.$

Put $\alpha>\beta\geq 0$ 
and  $I_n(y)=yI_{n-1}^\prime(y)$ ($n=1,2,3,\cdots$) and  $I_0^\prime(y)=I(y).$ 
Then 
$$I_n(y)=-y^\beta \beta^{n+1}(1+\beta)+e^{y-\alpha}y[y^{n+\beta}+\alpha^{1+\beta}\beta P_{n-1}(y)+y^\beta Q_{n-1}(y)]$$
with $P_{n-1}(y)=1+\sum\limits_{i=1}^{n-1}p_i^{n-1}y^i$ and 
$Q_{n-1}(y)=(1+\beta)^n+\sum\limits_{j=1}^{n-1}q_{n-j}^{n-1}(\beta)y^i$, 
where $p_i^{n-1}$ are positive integers,  
$q_j^{n-1}(\beta)$ are polynomials of degree $j$ 
whith respect to $\beta^j$ whose coefficients are positive integers, where  
$n=1,2,3,\cdots,$ $i=1,2,\cdots,n-1,$  $j=1,2,\cdots,n-1.$
It can be shown that  $I_n(\alpha)>0$ ($n=1,2,3,\cdots$) and that 
$$I_n(y)>-y^\beta\beta^{n+1}(1+\beta)+y^{\beta+1}(1+\beta)^n>y^\beta\beta(1+\beta)[(1+\beta)^{n-1}-\beta^n]$$
for all $y\geq \alpha>\beta\geq 0.$ 
Hence there exists a positive integer $n_0$ which depends only on $\beta,$ 
such that $I_n(y)>0$ for all $y\geq\alpha>\beta\geq 0$ as $n\geq n_0.$

By the definition of $I_n(y),$ we know that $I_{n_0-1}^\prime(y)>0$ for all 
$y\geq \alpha$ when $\alpha>\beta\geq 0.$ 
It follows that $I_{n_0-1}(y)>0$ for all 
$y\geq \alpha$ when $\alpha>\beta\geq 0.$ 
Repeating the above analysis, we conclude that $I_{n}^\prime(y)>0$ and $I_n(y)>I_n(\alpha)$ $(n=0,1,2,3,\cdots)$ for 
all $y\geq \alpha$ when $\alpha>\beta\geq 0.$  In particular, 
$I(y)>0$ for all $y\geq \alpha$ when $\alpha>\beta\geq 0.$

{\small 
\section*{Acknowledgement}
This paper is completed under the direction of Professor Zhu Xi-Ping. The authors would like to express their gratitude
for his continuous guidance and much valuable advice. XF is also grateful for support from NSFC 10171114. 
ZJ is supported by NSFC 10271121 and sponsored by SRF for ROCS, SEM.  
The authors would also like to thank Dr Chen Binglong   
for his helpful comments on this paper.  
Finally, the authors would like to thank the referee of this paper for his helpful comments and suggestions.

}

\end{document}